\documentclass[a4paper,11pt,makeidx]{amsart} 
\usepackage{amscd}
\usepackage{xypic}  
\usepackage{amssymb}
\usepackage{amsthm}
\usepackage{epsf}
\makeindex
\newtheoremstyle{break}
  {9pt}
  {9pt}
  {\itshape}
  {}
  {\bfseries}
  {.}
  {\newline}
  {}
\theoremstyle{break}
\newtheorem*{unico}{Theorem}

\theoremstyle{plain}
\newtheorem{thm}{Theorem}[section]
\newtheorem{cor}[thm]{Corollary}

\newtheorem{lemma}[thm]{Lemma}

\newtheorem{defn}[thm]{Definition}

\newtheorem{rem}[thm]{Remark}

\renewcommand{\proofname}{Proof}

\def\c1{\operatorname{c_1}}
\def\c2{\operatorname{c_2}}

\def\PP{{\mathbb P}}

\def\A{{\mathcal A}}

\def\L{{\mathcal L}}

\def\O{{\mathcal O}}

\def\*{\otimes}

\def\eqv{\equiv}

\def\+{\oplus}                   
\def\*{\otimes}                  

\def\Num{\operatorname{Num}}

\begin{document}

\title[A sharp vanishing theorem for line bundles on K3 or Enriques surfaces]{A sharp vanishing theorem
for line bundles \\ on K3 or Enriques surfaces} 

\author[A.L. Knutsen and A.F. Lopez]{Andreas Leopold Knutsen* and Angelo Felice Lopez**}

\address{\hskip -.43cm Andreas Leopold Knutsen, Dipartimento di Matematica, Universit\`a di Roma Tre, Largo San Leonardo
Murialdo 1, 00146, Roma, Italy. e-mail {\tt knutsen@mat.uniroma3.it}}

\thanks{* Research partially supported by a Marie Curie Intra-European Fellowship within the 6th European
Community Framework Programme}

\address{\hskip -.43cm Angelo Felice Lopez, Dipartimento di Matematica, Universit\`a di Roma 
Tre, Largo San Leonardo Murialdo 1, 00146, Roma, Italy. e-mail {\tt lopez@mat.uniroma3.it}}

\thanks{** Research partially supported by the MIUR national project ``Geometria delle variet\`a algebriche"
COFIN 2002-2004.}


\begin{abstract}
Let $L$ be a line bundle on a K3 or Enriques surface. We give a vanishing theorem for $H^1(L)$ that,
unlike most vanishing theorems, gives necessary and sufficient geometrical conditions for the vanishing.
This result is essential in our study of Brill-Noether theory of curves on Enriques surfaces \cite{kl1}
and of Enriques-Fano threefolds \cite{klm}.
\end{abstract}

\maketitle

\section{Introduction}
\label{intro}

\noindent Since Grothendieck's introduction of basic tools such as the cohomology of sheaves and the
Grothendieck-Riemann-Roch theorem, vanishing theorems have proved to be essential in many studies in
algebraic geometry.

\noindent Perhaps the most influential one, at least for line bundles, is the well-known
Kawamata-Viehweg vanishing theorem (\cite{ka, vie}) which, in its simplest form, asserts that
$H^i(K_X + \L) = 0$ for $i > 0$ and any big and nef line bundle $\L$ on a smooth variety $X$. On the
other hand, as most vanishing theorems (even for special surfaces \cite[Thm.1.5.1]{cd}), it gives only
{\it sufficient conditions for the vanishing}. Practice shows though that, in many situations, it would be
very useful to know that a certain vanishing is {\it equivalent} to some geometrical/numerical properties
of $\L$. 

\noindent In this short note we accomplish the above goal for line bundles on a K3 or
Enriques surface, by proving that, when $L^2 > 0$, the vanishing of $H^1(L)$ is equivalent to the fact
that the intersection of $L$ with all effective divisors of self-intersection $-2$ is at least $-1$. 

\noindent In the statement of the theorem we will employ the following

\begin{defn}
Let $X$ be a smooth surface. We will denote by $\sim$ (respectively $\eqv$) the linear (respectively
numerical) equivalence of divisors (or line bundles) on $X$. We will say that a line bundle
$L$ is {\bf primitive} if $L \eqv kL'$ for some line bundle $L'$ and some integer $k$ implies $k = \pm 
1$.
\end{defn} 

\begin{unico} 
Let $X$ be a K3 or an Enriques surface and let $L$ be a line bundle on $X$ such that 
$L > 0$ and $L^ 2 \geq 0$. Then $H^1(L) \neq 0$ if and only if one of the three
following occurs:
\begin{itemize}
\item[(i)] \ $L \sim nE$ for $E > 0$ nef and primitive with $E^2 = 0$, $n \geq 2$ and $h^1(L) = n - 1$ if
$X$ is a K3 surface, $h^1(L) = \lfloor \frac{n}{2} \rfloor$ if $X$ is an Enriques surface;
\item[(ii)] \ $L \sim nE + K_X$ for $E > 0$ nef and primitive with $E^2=0$, $X$ is an Enriques
surface, $n \geq 3$ and $h^1(L) = \lfloor \frac{n-1}{2} \rfloor$;
\item[(iii)] \ there is a divisor $\Delta > 0$ such that $\Delta^2 = -2$ and $\Delta.L \leq -2$.
\end{itemize}
\end{unico}

\noindent Note that the hypothesis $L > 0$ is not restrictive since, if $L$ is nontrivial, from $L^ 2 \geq
0$ we get by Riemann-Roch that either $L > 0$ or $K_X - L > 0$, and $h^1(L) = h^1(K_X - L)$ by Serre duality.

\noindent The theorem has of course many possible applications. For example, if $L$ is base-point
free and $|P|$ is an elliptic pencil on $X$, the knowledge of $h^0(L - nP)$ for $n \geq 1$ (which follows by
Riemann-Roch if we know that $h^1(L - nP) = 0$) determines the type of scroll spanned by the divisors of
$|P|$ in $\PP H^0(L)$ and containing $\varphi_L(X)$ (\cite{sd, kj, co}). Most importantly for
us, this result proves crucial in our study of the Brill-Noether theory \cite{kl1, kl2} and Gaussian
maps \cite{klGM} of curves lying on an Enriques surface, and especially in our proof of a genus bound for
threefolds having an Enriques surface as a hyperplane section given in \cite{klm}.

\vskip .2cm

\noindent {\it Acknowledgments}. The authors wish to thank Roberto Mu\~{n}oz for several helpful
discussions.

\section{Proof of the Theorem}

\noindent We first record the following simple but useful fact.

\begin{lemma} 
\label{lemma10}
Let $X$ be a smooth surface and let $A > 0$ and $B > 0$ be divisors on $X$ such that
$A^2 \geq 0$ and $B^2 \geq 0$. Then $A.B \geq 0$ with equality if and only if there exists a
primitive divisor $F > 0$ and integers $a \geq 1, b \geq 1$ such that $F^2 = 0$ and $A
\eqv aF, B \eqv bF$.
\end{lemma}

\begin{proof}
The first assertion follows from the signature theorem \cite[VIII.1]{bpv}. If $A.B = 0$, then we cannot
have $A^2 > 0$, otherwise the Hodge index theorem implies the contradiction $B \eqv 0$. Therefore $A^2 = B^2 =
0$. Now let $H$ be an ample line bundle on $X$ and set $\alpha = A.H, \beta = B.H$. We have $(\beta A -
\alpha B)^2 = 0$ and $(\beta A - \alpha B).H = 0$, therefore $\beta A \eqv \alpha B$ by the Hodge index
theorem. As there is no torsion in $\Num(X)$ we can find a divisor $F$ as claimed.
\end{proof}

\noindent We now proceed with the theorem.

\begin{proof}
One immediately sees that $h^1(L)$ has the given values in (i) and (ii). 
In the case (iii) we first observe that $h^2(L - \Delta) = 0$. In fact
$(K_X - L + \Delta)^2 > 0$, whence if $K_X - L + \Delta \geq 0$ the signature theorem
\cite[VIII.1]{bpv} implies $0 \leq L.(K_X - L + \Delta) = - L^2 + L.\Delta \leq -2$, a
contradiction. Therefore by Riemann-Roch we get
\[ \frac{1}{2}L^2 + \chi(\O_X) < \frac{1}{2}L^2 - \Delta.L -1 + \chi(\O_X) \leq h^0(L - \Delta) \leq h^0(L) = 
\frac{1}{2}L^2 + \chi(\O_X) + h^1(L) \]
whence $h^1(L) > 0$.

\noindent Now assume that $h^1(L) > 0$.

\noindent First we suppose that $L$ is nef. By Riemann-Roch we have that $L + K_X > 0$. Since $h^1(-
(L + K_X)) = h^1(L) > 0$, by \cite[Lemma12.2]{bpv}, we deduce that $L + K_X$ is not 1-connected,
whence that there exist $L' > 0$ and $L'' > 0$ such that $L + K_X \sim L'+ L''$ and $L'.L'' \leq 0$. Now
$(L')^2 \geq (L')^2 + L'.L'' = L'.L \geq 0$ and similarly $(L'')^2 \geq 0$, whence Lemma
\ref{lemma10} implies that $L' \eqv aE$, $L'' \eqv bE$ for some $a, b \geq 1$ and for $E > 0$ nef and
primitive with $E^2 = 0$. This gives us the two cases (i) and (ii).

\noindent Now assume that $L$ is not nef, so that the set
\[ \A_1(L):= \{ \Delta > 0 \; : \; \Delta^2 = -2, \  \Delta.L \leq -1 \} \]
is not empty. Similarly define the set
\[ \A_2(L) = \{ \Delta > 0 \; : \; \Delta^2 = -2, \ \Delta.L \leq -2 \}. \]

\noindent If $\A_2(L) \neq \emptyset$ we are done. Assume therefore that $\A_2(L) = \emptyset$
and pick $\Gamma \in \A_1(L)$. Then $\Gamma.L = -1$, and we can clearly assume that $\Gamma$ is
irreducible. Hence if we set $L_1 = L- \Gamma$ we have that $L_1 > 0$, $L_1^2 = L^2$ and, since
$h^0(L_1) = h^0(L)$, also that $h^1(L_1) = h^1(L) > 0$.

\noindent If $L_1$ is nef, by what we have just seen, we have $L_1 \eqv nE$, for $n \geq
2$, whence $L \eqv nE + \Gamma$ and $-1 = \Gamma.L = nE.\Gamma - 2$, a contradiction.

\noindent Therefore $L_1$ is not nef and $\A_1(L_1) \neq \emptyset$.

\noindent If $\A_2(L_1) \neq \emptyset$ we pick a $\Delta \in \A_2(L_1)$. We have 
$-2 \geq \Delta.L_1=\Delta.(L-\Gamma) \geq -1-\Delta.\Gamma$, whence $\Delta.\Gamma \geq 1$,
$(\Delta+\Gamma)^2 \geq -2$ and $(\Delta+\Gamma).L_1 \leq -1$. Now Lemma \ref{lemma10}
yields $(\Delta+\Gamma)^2 = -2$, so that $\Delta.\Gamma =1$. Also $-1 \leq \Delta.L =
\Delta.(L_1 + \Gamma) \leq -1$, whence $\Delta.L = -1$ and $(\Delta + \Gamma).L = -2$, contradicting 
$\A_2(L) = \emptyset$.

\noindent We have therefore shown that $\A_2(L_1) = \emptyset$.

\noindent This means that we can continue the process. But the process must eventually stop,
since we always remove base components. This gives the desired contradiction. 
\end{proof}
\renewcommand{\proofname}{Proof}

\begin{rem} A naive guess, to insure the vanishing of $H^1(L)$ for a line bundle $L > 0$ with
$L^ 2 \geq 0$, could be that it is enough to add the hypothesis $L.R \geq -1$ for every irreducible rational
curve $R$. However this is not true. {\rm Take, for example, a nef divisor $B$ with $B^2 \geq 4$
and two irreducible rational curves $R_1, R_2$ such that $B.R_i = 0, R_1.R_2 = 1$. Then $L:= B + R_1 + R_2$
satisfies the above requirements, but $L.(R_1 + R_2)=-2$, whence $H^1(L) \neq 0$ by the theorem.}
\end{rem}

\begin{rem} {\rm It would be of interest to know if, in the statement of the theorem, it is
possible to replace divisors $\Delta > 0$ such that $\Delta^2 = -2$ with chains of irreducible rational
curves.}
\end{rem}

\begin{defn} 
\label{def:qnef}
An effective line bundle $L$ on a K3 or Enriques surface is said to be {\bf quasi-nef} if $L^2
\geq 0$ and $L.\Delta \geq -1$ for every $\Delta$ such that $\Delta > 0$ and $\Delta^2 = -2$. 
\end{defn}

\noindent An immediate consequence of the theorem is

\begin{cor} 
\label{cor:qnef} 
An effective line bundle $L$ on a K3 or Enriques surface is quasi-nef if and only if $L^2 \geq 0$
and either $h^1(L) = 0$ or $L \eqv nE$ for some $n \geq 2$ and some primitive and nef divisor $E > 0$
with $E^2 = 0$.
\end{cor}

\end{document}